\documentclass[12pt]{article}
\usepackage{amssymb,amsfonts,amsmath, psfrag,eepic,colordvi,epsfig}
by \topskip 
\numberwithin{equation}{section}
\newtheorem{theorem}{Theorem}[section]

\newtheorem{corollary}[theorem]{Corollary}

\newtheorem{property}[theorem]{Property}

\begin{document}
\parskip 6pt

\pagenumbering{arabic}
\def\sof{\hfill\rule{2mm}{2mm}}
\def\ls{\leq}
\def\gs{\geq}
\def\SS{\mathcal S}
\def\qq{{\bold q}}
\def\MM{\mathcal M}
\def\TT{\mathcal T}
\def\EE{\mathcal E}
\def\lsp{\mbox{lsp}}
\def\rsp{\mbox{rsp}}
\def\pf{\noindent {\it Proof.} }
\def\mp{\mbox{pyramid}}
\def\mb{\mbox{block}}
\def\mc{\mbox{cross}}
\def\qed{\hfill \rule{4pt}{7pt}}
\def\block{\hfill \rule{5pt}{5pt}}
\begin{center}
{\Large {\bf The Butterfly Decomposition of Plane Trees}}

\vskip 6mm

{\small William Y.C. Chen$^1$, Nelson Y. Li$^2$, and Louis W. Shapiro$^3$\\[%
2mm] $^{1,2}$Center for Combinatorics, LPMC\\
Nankai University,
Tianjin 300071,
P.R. China \\[3mm]
$^3$Howard University, Washington, DC 20059, USA \\[3mm]
$^1$chen@nankai.edu.cn, $^2$nelsonli@eyou.com, $^3$lshapiro@howard.edu \\[0pt%
]
}
\end{center}

\vskip10mm
\noindent {\bf Abstract.} We introduce the notion of doubly rooted
plane trees and give a decomposition of these trees, called the
butterfly decomposition which turns out to have many applications.
From the butterfly decomposition we obtain a one-to-one
correspondence between doubly rooted plane trees and free Dyck
paths, which implies a simple derivation of a relation between the
Catalan numbers and the central binomial coefficients. We also
establish a one-to-one correspondence between leaf-colored doubly
rooted plane trees and free Schr\"oder paths.  The classical
Chung-Feller theorem on free Dyck paths and some generalizations
and variations with respect to Dyck paths and Schr\"oder paths
with flaws turn out to be immediate consequences of the butterfly
decomposition and the preorder traversal of plane trees. We obtain
two involutions on free Dyck paths and free Schr\"oder paths,
leading to two combinatorial identities. We also use the butterfly
decomposition to give a combinatorial treatment of the generating
function for the number of chains in plane trees due to Klazar. We
further study the average size of chains in plane trees with $n$
edges and show that this number asymptotically tends to ${ n+9
\over 6}$.




\vskip10mm

\section{Introduction}

This paper is concerned with the enumeration of plane trees and
the number of chains in plane trees with $n$ edges. Although this
subject has been very well studied over many decades, it seems
that interesting problems still emerge. As we shall see, the
enumeration of chains in plane trees leads us to discover a
fundamental property of doubly rooted plane trees which has many
applications. We call this the butterfly decomposition.

From the butterfly decomposition, we can establish a
correspondence between doubly rooted plane trees and free Dyck
paths. So we immediately get the relation between the Catalan
numbers and the central binomial coefficients. The butterfly
decomposition also implies the classical Chung-Feller theorem  on
free Dyck paths with a given number of steps under $x$-axis. The
Chung-Feller theorem was first proved by Major Percy A. MacMahon
in 1909 \cite[p.168]{mac} but named after its 1949 re-discoverers
\cite{chung}. MacMahon proved it using formal series of words on
an alphabet; Chung and Feller used generating functions. The
previous combinatorial approaches to the Chung-Feller theorem are
based on the cycle lemma or cyclic paths, see Dershowitz-Zaks
\cite{derzak} and Narayana \cite{narayana}. There are other
Chung-Feller type results and generalizations in
\cite{callan,cameron,efy,eu,shapwoan,woan}. In some sense, the
butterfly decomposition can be regarded as labelled variation of
the cycle lemma.

The butterfly decomposition  also leads to the following results:
a correspondence between leaf-colored doubly rooted plane trees
and free Schr\"oder paths, a simple bijection between leaf-colored
plane trees and Schr\"oder paths, and a combinatorial
interpretation of the generating function for the number of chains
in plane trees obtained by Klazar \cite{klazar}. We show that
there is a one-to-one correspondence between chains in plane trees
and tricolored plane trees (the definition is given in Section 5).

We obtain two involutions on free Dyck paths and free Schr\"oder
paths, which lead to two combinatorial identities. The last
section of this paper gives an asymptotic formula for the average
size of chains in plane trees with $n$ edges.

\section{The Butterfly Decomposition}

In this section, we introduce the notion of doubly rooted plane
trees and their butterfly decomposition. This decomposition seems
to be fundamental for the enumeration of plane trees. It also
implies the generating function for the number of chains in plane
trees obtained by Klazar \cite{klazar}. We will study the
enumeration of chains in Section 5. The main result of this
section is a correspondence between doubly rooted plane trees and
free Dyck paths, from which it follows a combinatorial
interpretation of the relation
\begin{equation} \label{c-n-d}
(n+1) c_n  = {2n \choose n}.
\end{equation}
We will also establish a correspondence between free Dyck paths
and $2$-colored plane trees.

\begin{figure}[h,t]
\begin{center}
\setlength{\unitlength}{1mm}
\begin{picture}(25,15)

\multiput(10,10)(0,-10){2}{\circle*{1}}
\put(10,10){\line(-2,-1){10}} \put(10,10){\line(-1,-2){5}}
\put(0,5){\line(1,-1){5}} \put(10,10){\line(0,-1){10}}
\put(3.5,4){\footnotesize{$L_i$}}
\put(9.5,12){\footnotesize{$v_i$}} \put(10,10){\line(2,-1){10}}
\put(10,10){\line(1,-2){5}} \put(20,5){\line(-1,-1){5}}
\put(13.5,4){\footnotesize{$R_i$}}

\end{picture}
\end{center}
\caption{A Butterfly} \label{butterfly}
\end{figure}
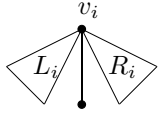

A (rooted) plane tree $T$ with a distinguished vertex $w$ is
called a {\em doubly rooted plane tree}, where the distinguished
vertex is regarded as the second root. The {\em butterfly
decomposition} of a doubly rooted plane tree $T$ with a
distinguished vertex $w$ is described as follows. Let $P= v_1 v_2
\ldots, v_k w$ be the path from the root of $T$ to $w$. Let $L_1,
L_2, \ldots, L_k$ be the subtrees such that $L_i$ consists of the
vertex $v_i$ and its descendants on the left hand side of the path
$P$. Similarly, we can define the subtrees $R_1, R_2, \ldots, R_k$
as the subtrees rooted at $v_1, v_2, \ldots, v_k$ consisting of
the descendants on the right hand side of $P$. Moreover, the
subtree of $T$ rooted at $w$ is denoted by $T'$. Therefore, a
plane tree $T$ with a distinguished vertex $w$ can be decomposed
into smaller structures  $(U_1, U_2, \ldots, U_k; T')$, where
$U_i$ is called a {\em butterfly} consisting of $L_i$ and $R_i$
and the edge in the middle, as shown in Figure~\ref{butterfly}.

Let $B$ and $C$ be the generating functions of the central
binomial coefficients and Catalan numbers respectively:
\begin{equation} \label{db}
B=\frac{1}{\sqrt{1-4x}}=\sum_{n=0}^{\infty }{\binom{2n}{n}}x^{n},
\end{equation}
\begin{equation}\label{dc}
C=\frac{1-\sqrt{1-4x}}{2x}=\sum_{n=0}^{\infty
}\frac{1}{n+1}{\binom{2n}{n}} x^{n},
\end{equation}
where ${2n \choose n}$ is called the {\em central binomial
coefficient}, and $c_n={1 \over n+1} {2n \choose n}$ is the $n$-th
{\em Catalan number} \cite{stanley}.
 It is clear that the generating function for a
butterfly with $n$ edges equals $xC^2$, and the generating
function for a sequence of $k$ butterflies with a total number of
$n$ edges equals $(xC^2)^k$. Note that the number of doubly rooted
plane trees with $n$ edges equals  $n+1$ times the Catalan number,
that is, the central binomial coefficient ${2n \choose n}$. Thus,
we arrive at the following generating function relation:
\begin{equation}\label{g-b}
B= C+ C(xC^2) + C(xC^2)^2 + \cdots = {C\over 1- xC^2}.
\end{equation}

A natural question arises: is there a simple combinatorial
argument that leads to this conclusion without resorting to the
formula for the Catalan numbers? The answer is affirmative, this
leads to a quite simple derivation of the relation (\ref{c-n-d}).

Recall that a {\em Dyck path} of length $2n$ is a lattice path
from the origin to $(2n,0)$ consisting of up steps $(1, 1)$ and
down steps $(1, -1)$ that does not go below the $x$-axis. An {\em
elevated Dyck path} or {an irreducible Dyck path} is defined as a
Dyck path that does not touch the $x$-axis except for the origin
and the destination. A lattice path from the origin to $(2n,0)$
using the steps $(1, 1)$ and $(1, -1)$ without additional
restrictions is called \emph{a free Dyck path}, a free Dyck path
is also called a Dyck path with flaws in the sense that the
segments below the $x$-axis are regarded as flaws, see  Eu-Liu-Yeh
\cite{eu}. The reflection of a Dyck path with respect to the
$x$-axis is called a negative Dyck path. An {\em elevated
(irreducible) negative Dyck path} is defined in the same manner.
As we shall see, free Dyck paths can be regarded as a labelled
version of Dyck paths. Clearly, the set of free Dyck paths of
length $2n$ is just the set of sequences consisting of $n$ up
steps and $n$ down steps, as counted by the central binomial
coefficient ${2n \choose n}$.

\begin{theorem} \label{t-c-n}
There is a bijection between the set of doubly rooted plane trees
with $n$ edges and the set of free Dyck paths of length $2n$.
\end{theorem}

First we give a combinatorial setting for the proof of the above
theorem. We recall the classical glove bijection between plane
trees and Dyck paths \cite {dershowitz}. This correspondence is
also referred to as the preorder traversal of a plane tree. For
the purpose of this paper, we may view the glove bijection as a
recursive procedure. Recall that a {\em planted plane tree} is a
plane tree whose root has only one child. Then the glove bijection
gives a correspondence between the set of planted plane trees with
$n$ edges and the set of elevated Dyck paths of length $2n$. A
planted plane tree with one edge corresponds to the elevated Dyck
path of length two. Let $T$ be a planted plane tree, and let $T_1,
T_2, \ldots, T_k$ be the  subtrees of the only one child of the
root of $T$. Let $P_1, P_2, \ldots, P_k$ be elevated Dyck paths
corresponding to $T_1, T_2, \ldots, T_k$ respectively. Then
$UP_1P_2\cdots P_kD$ is an elevated Dyck path of length $2n$,
where $U$ stands for an up step and $D$ stands for a down step.

We are now ready to give a proof of Theorem \ref{t-c-n}.

\pf Let $T$ be a doubly rooted plane tree with $n$ edges. Let $w$
be the distinguished vertex of $T$ and let $v_1v_2\cdots v_k w$ be
the path from the root to $w$. Suppose that $(L_1, R_1; L_2, R_2;
\ldots ; L_k, R_k; T')$ is the butterfly decomposition of $T$.

For the $L_i$ $(1 \leq i \leq k)$ and $T'$, we use the glove
bijection to them and call the resulting Dyck paths $P_i$ and
$P_{k+1}$. For every $R_i$, we first add an edge at each root to
form a planted plane tree $T_i$, then use the glove bijection to
produce a negative elevated Dyck path $Q_i$. Now
 \begin{equation} \label{fdf}
  P_1 \, Q_1 \, P_2\, Q_2\, \cdots \, P_k \, Q_k \, P_{k+1}
  \end{equation}
 is a free Dyck path of length $2n$.  Conversely, given a free
 Dyck path we may decompose it into elevated (irreducible)
 segments  like the {\em first return decomposition}
   of a
Dyck path \cite{deutschd}, and we may reverse the above procedure
to construct a doubly rooted plane tree because any free Dyck path
$P$ has a unique decomposition in the form (\ref{fdf}) such that
$Q_1, Q_2, \ldots, Q_k$ are negative elevated Dyck paths and $P_1,
P_2, \ldots, P_{k+1}$ are the usual Dyck paths with the empty
paths allowed. Thus we have established the bijection. \qed

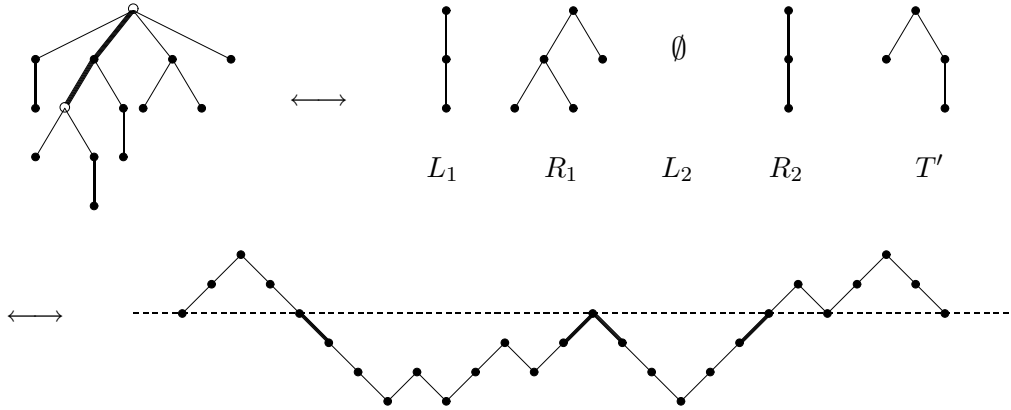
\begin{figure}[h,t]
\begin{center}
\setlength{\unitlength}{1.3mm}
\begin{picture}(90,45)

\put(10,40.2){\circle{1}} \put(10,40){\line(-2,-1){10}}
\put(10,40){\line(-4,-5){4}}\put(9.7,40){\line(-4,-5){4}}
\put(9.8,40){\line(-4,-5){4}}\put(10.1,40){\line(-4,-5){4}}
\put(9.9,40){\line(-4,-5){4}}\put(10.2,40){\line(-4,-5){4}}
\put(10,40){\line(4,-5){4}} \put(10,40){\line(2,-1){10}}
\multiput(0,35)(6,0){2}{\circle*{.8}}
\multiput(14,35)(6,0){2}{\circle*{.8}}

\put(0,35){\line(0,-1){5}} \put(0,30){\circle*{.8}}

\put(6,35){\line(-3,-5){2.7}}\put(5.7,35){\line(-3,-5){2.7}}
\put(5.8,35){\line(-3,-5){2.7}}\put(6.1,35){\line(-3,-5){2.8}}
\put(5.9,35){\line(-3,-5){2.7}}\put(6.2,35){\line(-3,-5){2.8}}
\put(6,35){\line(3,-5){3}} \put(3,30.1){\circle{1}}
\put(9,30){\circle*{.8}}

\put(14,35){\line(-3,-5){3}} \put(14,35){\line(3,-5){3}}
\put(11,30){\circle*{.8}} \put(17,30){\circle*{.8}}

\put(3,30){\line(-3,-5){3}} \put(3,30){\line(3,-5){3}}
\put(0,25){\circle*{.8}} \put(6,25){\circle*{.8}}

\put(9,30){\line(0,-1){5}} \put(9,25){\circle*{.8}}

\put(6,25){\line(0,-1){5}} \put(6,20){\circle*{.8}}

\put(26,30){$\longleftrightarrow$}


\put(42,40){\line(0,-1){10}}
\multiput(42,40)(0,-5){3}{\circle*{.8}} \put(40,23){\small{$L_1$}}

\put(55,40){\circle*{.8}} \put(55,40){\line(-3,-5){3}}
\put(55,40){\line(3,-5){3}} \multiput(52,35)(6,0){2}{\circle*{.8}}
\put(52,35){\line(-3,-5){3}} \put(52,35){\line(3,-5){3}}
\multiput(49,30)(6,0){2}{\circle*{.8}}\put(52,23){\small{$R_1$}}

\put(65,35){$\emptyset$} \put(64,23){\small{$L_2$}}

\put(77,40){\line(0,-1){10}}
\multiput(77,40)(0,-5){3}{\circle*{.8}} \put(75,23){\small{$R_2$}}

\put(90,40){\circle*{.8}} \put(90,40){\line(-3,-5){3}}
\put(90,40){\line(3,-5){3}} \multiput(87,35)(6,0){2}{\circle*{.8}}
\put(93,35){\line(0,-1){5}}
\put(93,30){\circle*{.8}}\put(90,23){\small{$T'$}}

\put(-3,7.9){$\longleftrightarrow$}

\multiput(10,9)(1,0){90}{\line(1,0){0.5}}

\put(15,9){\line(1,1){6}}\put(21,15){\line(1,-1){15}}
\put(36,0){\line(1,1){3}} \put(39,3){\line(1,-1){3}}
\put(42,0){\line(1,1){6}} \put(48,6){\line(1,-1){3}}
\put(51,3){\line(1,1){6}} \put(57,9){\line(1,-1){9}}
\put(66,0){\line(1,1){12}}\put(78,12){\line(1,-1){3}}
\put(81,9){\line(1,1){6}}\put(87,15){\line(1,-1){6}}

\put(15,9){\circle*{.8}}\put(18,12){\circle*{.8}}\put(21,15){\circle*{.8}}
\put(24,12){\circle*{.8}}\put(27,9){\circle*{.8}}\put(30,6){\circle*{.8}}
\put(33,3){\circle*{.8}}\put(36,0){\circle*{.8}}\put(39,3){\circle*{.8}}
\put(42,0){\circle*{.8}}\put(45,3){\circle*{.8}}\put(48,6){\circle*{.8}}
\put(51,3){\circle*{.8}}\put(54,6){\circle*{.8}}\put(57,9){\circle*{.8}}
\put(60,6){\circle*{.8}}\put(63,3){\circle*{.8}}\put(66,0){\circle*{.8}}
\put(69,3){\circle*{.8}}\put(72,6){\circle*{.8}}\put(75,9){\circle*{.8}}
\put(78,12){\circle*{.8}} \put(81,9){\circle*{.8}}
\put(84,12){\circle*{.8}}\put(87,15){\circle*{.8}}
\put(90,12){\circle*{.8}}\put(93,9){\circle*{.8}}

\put(27.3,9){\line(1,-1){3}}\put(54.3,6){\line(1,1){3}}
\put(57.3,9){\line(1,-1){3}}\put(72.3,6){\line(1,1){3}}
\put(27.2,9){\line(1,-1){3}}\put(54.2,6){\line(1,1){3}}
\put(57.2,9){\line(1,-1){3}}\put(72.2,6){\line(1,1){3}}
\put(27.1,9){\line(1,-1){3}}\put(54.1,6){\line(1,1){3}}
\put(57.1,9){\line(1,-1){3}}\put(72.1,6){\line(1,1){3}}
\put(26.9,9){\line(1,-1){3}}\put(53.9,6){\line(1,1){3}}
\put(56.9,9){\line(1,-1){3}}\put(71.9,6){\line(1,1){3}}
\end{picture}
\end{center}
\caption{Doubly rooted plane trees and free Dyck paths} \label{df}
\end{figure}

An example of the above bijection is shown in Figure~ \ref{df}.

We next give another interpretation of the generating function for
the number of bicolored plane trees.
 Guided by
the following generating function identity
\begin{equation}\label{i-1-2}
 { C \over 1 -x C^2} = {1 \over 1-2xC},
 \end{equation}
 we are led
to introduce the notion of bicolored plane trees and $k$-colored
plane trees, in general. A {\em $k$-colored plane tree} is a plane
tree in which the children of the root are colored with $k$
colors. A $2$-colored plane tree is called a {\em bicolored plane
tree}, and a $3$-colored plane tree is called a {\em tricolored
plane tree}. For bicolored plane trees, we assume that the two
colors are black and white.
 Note that this terminology is somewhat
misleading because in our context only the children of the root
are colored.  The relation (\ref{i-1-2}) indicates that the set of
bicolored plane trees are in one-to-one correspondence with doubly
rooted plane trees. We next establish such a correspondence by
making a connection between bicolored plane trees and free Dyck
paths.

\begin{theorem} \label{t-c-n-2}
There is a one-to-one correspondence between the set of bicolored
plane trees with $n$ edges and the set of free Dyck paths of
length $2n$.
\end{theorem}

\pf  Let $T$ be a bicolored plane tree, and let $T_1, T_2, \ldots,
T_k$ be the planted subtrees of the root of $T$, listed from left
to right. If $T_i$ inherits the black color, then we construct an
negative elevated Dyck path $P_i$ from $T_i$; otherwise we
construct an elevated Dyck path $P_i$ above the $x$-axis. So we
get a free Dyck path $P_1P_2\ldots P_k$. Conversely, given a free
Dyck path we may construct a bicolored plane tree. Hence we obtain
the desired bijection. \qed

The bijections in Theorems \ref{t-c-n} and \ref{t-c-n-2} lead to a
bijection between doubly rooted plane trees and bicolored plane
trees. In fact, we may establish a direct correspondence without
resorting to free Dyck paths.

\begin{theorem}\label{doubly}
There is a bijection between the set of doubly rooted plane trees
with $n$ edges and the set of bicolored plane trees with $n$
edges.
\end{theorem}

\pf By the butterfly decomposition in Theorem~\ref{t-c-n}, we get
subtrees $L_i$, $T_i$ and $T'$. By coloring $L_i$ and $T'$ black
while $T_i$ white, and identifying their roots as the root of the
corresponding bicolored plane tree, we have its subtrees listed
from left to right as
\[ L_1 \, T_1 \, L_2\, T_2\, \cdots \, L_k \, T_k \, T'.\]
The reverse procedure is easy to construct. Thus we have
established the bijection. \qed

\section{The Chung-Feller Theorem}

We begin this section by pointing out that the classical
Chung-Feller theorem on Dyck paths is an immediate consequence of
our bijection between doubly rooted plane trees and free Dyck
paths. To see this connection, one only needs a simple observation
on the preorder traversal of a plane tree. We also use this idea
to derive some refinements and generalizations of the Chung-Feller
theorem, including some recent results of Eu, Fu and Yeh
\cite{efy} on Dyck paths and Schro\"oder paths with flaws.

\begin{theorem}[Chung-Feller]\label{cf}
For any $0\leq m \leq n$, the number of free Dyck paths of length
$2n$ that contain exactly $2m$ steps below the $x$-axis is
independent of $m$, and is equal to the $n$-th Catalan number
$c_n$.
\end{theorem}

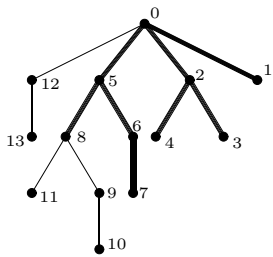
\begin{figure}[h,t]
\begin{center}
\setlength{\unitlength}{1.5mm}
\begin{picture}(20,22)

\put(10,20){\circle*{.8}}\put(10,20){\line(-2,-1){10}}
\put(10.5,20.5){\tiny{$0$}}\put(20.5,15.5){\tiny{$1$}}
\put(9.8,20){\line(-4,-5){4}}\put(9.9,20){\line(-4,-5){4}}
\put(10,20){\line(-4,-5){4}}\put(10.1,20){\line(-4,-5){4}}
\put(10.2,20){\line(-4,-5){4}}
\put(9.8,20){\line(4,-5){4}}\put(9.9,20){\line(4,-5){4}}
\put(10,20){\line(4,-5){4}}\put(10.1,20){\line(4,-5){4}}
\put(10.2,20){\line(4,-5){4}}
\put(9.7,20){\line(2,-1){10}}\put(9.8,20){\line(2,-1){10}}
\put(9.9,20){\line(2,-1){10}}\put(10,20){\line(2,-1){10}}
\put(10.1,20){\line(2,-1){10}}\put(10.2,20){\line(2,-1){10}}
\put(10.3,20){\line(2,-1){10}}\put(10.4,20){\line(2,-1){10}}
\multiput(0,15)(6,0){2}{\circle*{.8}}
\multiput(14,15)(6,0){2}{\circle*{.8}}

\put(0,15){\line(0,-1){5}} \put(0,10){\circle*{.8}}

\put(6.8,14.5){\tiny{$5$}}
\put(5.8,15){\line(-3,-5){3}}\put(5.9,15){\line(-3,-5){3}}
\put(6,15){\line(-3,-5){3}}\put(6.1,15){\line(-3,-5){3}}
\put(6.2,15){\line(-3,-5){3}}
\put(5.8,15){\line(3,-5){3}}
\put(5.9,15){\line(3,-5){3}}\put(6,15){\line(3,-5){3}}
\put(6.1,15){\line(3,-5){3}}\put(6.2,15){\line(3,-5){3}}
\put(3,10){\circle*{.8}} \put(9,10){\circle*{.8}}

\put(14.5,15){\tiny{$2$}}\put(11.8,9){\tiny{$4$}}\put(17.8,9){\tiny{$3$}}
\put(13.8,15){\line(-3,-5){3}}\put(13.9,15){\line(-3,-5){3}}
\put(14,15){\line(-3,-5){3}}\put(14.1,15){\line(-3,-5){3}}
\put(14.2,15){\line(-3,-5){3}}
\put(13.8,15){\line(3,-5){3}}\put(13.9,15){\line(3,-5){3}}
\put(14,15){\line(3,-5){3}}\put(14.1,15){\line(3,-5){3}}
\put(14.2,15){\line(3,-5){3}}
\put(11,10){\circle*{.8}} \put(17,10){\circle*{.8}}

\put(3,10){\line(-3,-5){3}} \put(3,10){\line(3,-5){3}}
\put(0,5){\circle*{.8}} \put(6,5){\circle*{.8}}

\put(8.9,10.5){\tiny{$6$}}\put(9.5,4.5){\tiny{$7$}}\put(4,9.5){\tiny{$8$}}
\put(8.8,10){\line(0,-1){5}}\put(8.9,10){\line(0,-1){5}}
\put(9,10){\line(0,-1){5}}\put(9.1,10){\line(0,-1){5}}\put(9.2,10){\line(0,-1){5}}
\put(9,5){\circle*{.8}}

\put(6,5){\line(0,-1){5}} \put(6,0){\circle*{.8}}

\put(6.7,4.5){\tiny{$9$}}\put(6.7,0){\tiny{$10$}}\put(.7,4.2){\tiny{$11$}}
\put(.8,14.2){\tiny{$12$}}\put(-2.3,9.2){\tiny{$13$}}
\end{picture}
\end{center}
\caption{Labels for the Chung-Feller theorem} \label{chung}
\end{figure}

Using the butterfly decomposition, we may transform the
Chung-Feller theorem to an equivalent form on plane trees, which
turns out to be a simple property of the preorder traversal. To be
precise, we define the {\em right-to-left preorder traversal} of a
plane tree $T$ as a recursive procedure. First, visit the root of
$T$. Let $T_1, T_2, \ldots, T_k$ be the subtrees of the root of
$T$ listed from left to right. Then traverse the subtrees in the
order of $T_k, T_{k-1}, \ldots, T_{1}$. From the above traversal
procedure, we may label the vertices of $T$ with the numbers $0,
1, 2, \ldots, n$ in the order that they are visited.
Figure~\ref{chung} gives the plane tree corresponding to the free
Dyck path in Figure~\ref{df} and the labels by the right-to-left
preorder traversal.

The following property immediately implies the Chung-Feller
theorem since any plane tree can be regarded as a doubly rooted
plane tree in which the distinguished vertex is chosen as the
vertex with a given label $m$ with respect to the right-to-left
preorder traversal.

\begin{theorem} Let $T$ be a plane tree with $n$ edges. Assume
that the vertices of $T$ are labelled by $0, 1, 2, \ldots, n$
according to the right-to-left preorder traversal. Let $w$ be the
vertex labelled by $m$, where $m$ is a given number not exceeding
$n$. Then the doubly rooted plane tree $T$ with $w$ being the
distinguished vertex corresponds to a free Dyck path with $m$ down
steps (up steps) below the $x$-axis.
\end{theorem}

As another corollary, we note that half of all free Dyck paths end
with an up step. Thus over all plane trees with $n$ edges, half of
the vertices are leafs, see Problem 10753 of the American
Mathematical Monthly \cite{Hunger, shapiro-p}.


The above interpretation of the Chung-Feller theorem also implies
some refinements and generalizations  recently obtained by Eu, Fu
and Yeh \cite{efy}.   Let us define some terminology. We say that
a free Dyck path has $m$ {\em flaws} if it contains $m$ up (or
down) steps below the $x$-axis. We note that a negative elevated
(irreducible) Dyck path is called a {\em flaw block} by Eu, Fu and
Yeh \cite{efy}. We define the {\em stem} of a doubly rooted plane
tree as the path from the root to the distinguished vertex. Let
$T$ be a doubly rooted plane tree with a distinguished vertex $w$.
An edge of $T$ is said to be a {\em prefix edge} if it is either
on on the stem of $T$ or to the right of the stem. In other words,
a prefix edge is an edge with labels not exceeding the label of
the distinguished vertex with respect to the right-to-left
preorder traversal. An example is shown in Figure~\ref{chung}
where the prefix edges are drawn with thick edges.

Using the preorder traversal of plane trees, we get the following
generalization of the refined version of the Chung-Feller theorem
\cite{efy}.

\begin{theorem}
For $0\leq k \leq m \leq n$, there is a bijection between the set
of free Dyck path of length $2n$ with $m$ flaws in $k$ flaw blocks
and the set of doubly rooted plane trees of $n$ edges with stem
size $k$ and $m$ prefix edges.
\end{theorem}

\pf From the butterfly decomposition and the correspondence in
Theorem~\ref{t-c-n}, we see that the number of flaw blocks in a
free Dyck path equals the stem size of the corresponding doubly
rooted plane tree, and the number of flaws in a free Dyck path
equals the number of prefix edges in the plane tree. This
completes the proof. \qed

By the butterfly decomposition, one sees that the generating
function for doubly rooted plane trees with stem size $k$ equals
$x^kC^k\cdot C^{k+1}$. It follows that the number of such trees
with $n$ edges and $m$ prefix edges equals
$[x^m]x^kC^k\cdot[x^{n-m}]C^{k+1}$, where $[x^n]C^k$ is the usual
notation for the coefficient of $x^n$ in the expansion of $C^k$.
By the Lagrange inversion formula \cite{stanley}, we have
\begin{equation} \label{xn}
[x^n]C^k=\frac{k}{2n+k}{2n+k \choose n}.
\end{equation}
 Thus, we obtain the following expression.

\begin{corollary}
For $0< k \leq m \leq n$, the number of free Dyck paths of length
$2n$ with $m$ flaws and $k$ flaw blocks equals
\[
\frac{k}{2m-k}{2m-k \choose m}{k+1 \over 2n-2m+k+1} {2n-2m+k+1
\choose n-m}.
\]
\end{corollary}

Setting $m=n$ in the above corollary, one gets the number of Dyck
paths of length $2n$ with $k$ returns obtained by Deutsch
\cite{deutschd}:
\begin{equation}\label{cf-i-2}
\frac{k}{2n-k}{2n-k \choose n}.
\end{equation}


We next consider the enumeration of Schr\"oder paths with flaws.
For this purpose, we need to introduce the notion of {\em
leaf-colored doubly rooted plane trees} which are defined as
doubly rooted plane trees whose leaves are colored with two colors
red $(R)$ and blue $(B)$ under the convention that the
distinguished vertex receives no color even if it is a leaf. An
edge of a plane tree is called an {\em external edge} if it's end
vertex is a leaf; otherwise it is called an {\em internal edge}.
As we shall see, such leaf-colored doubly rooted plane trees are
in one-to-one correspondence with free Schr\"oder paths. We note
that there is a bijection between Schr\"oder paths and plane trees
with every leaf being colored red or blue, see Gouyou-Beauchamps
and Vauquelin \cite{guo}, here we create a new one.

Recall that a {\em Schr\"oder path} of length $2n$ is a lattice
path in the plane from $(0,0)$ to $(2n,0)$ with up steps
$U=(1,1)$, horizontal steps $H=(2,0)$, and down steps $D=(1,-1)$,
that never go below the $x$-axis. These paths are enumerated by
the Schr\"oder numbers $r_n$ \cite{stanley}. An {\em elevated
(irreducible) Schr\"oder path} and a {\em negative Schr\"oder
path} are defined in the same manner as with Dyck paths. A lattice
path from $(0,0)$ to $(2n,0)$ with steps $U=(1,1)$, $H=(2,0)$, and
$D=(1,-1)$ without additional restrictions is called a {\em free
Schr\"oder path}. We say that a free Schr\"oder path has $m$ {\em
flaws} if the number of $U$ steps and  $H$ steps under the
$x$-axis equals $m$. A {\em flaw block} of a Schr\"oder path is
defined as  a negative elevated Schr\"oder path.

By the  preorder traversal, we obtain the following
correspondence.

\begin{theorem}
There is a one-to-one correspondence between the set of plane
trees with $n$ edges in which each leaf  is colored red or blue
and the set of Schr\"oder paths of length $2n$.
\end{theorem}

\pf  Let $T$ be a plane tree with $n$ edges in which each leaf is
colored red or blue. We proceed to construct a Schr\"oder path of
length $2n$ from the (left-to-right) preorder traversal. In the
preorder traversal of the vertices of $T$, each edge is visited
twice. Note that when an external edge $e=(u,v)$ ($v$ is a leaf)
is traversed, one always visits the vertex $u$, then the leaf $v$,
and then immediately goes back to the vertex $u$. Now we may
generate a sequence of $U$, $D$, and $H$ steps by the following
rule: (1) When an internal edge is visited for the first time, we
get an $U$ step. (2) When an internal edge is visited for the
second time, we get a $D$ step. (3) When an external edge with a
red leaf is traversed, we get two steps $UD$. (4) When an external
edge with a blue leaf is traversed, we get an $H$ step. It is easy
to see that we obtain a Schr\"oder path of length $2n$ and the
above procedure is reversible. \qed

By using the butterfly decomposition, we obtain the following
correspondence.

\begin{theorem} \label{t-s-n}
There is a bijection between the set of leaf-colored doubly rooted
plane trees with $n$ edges and the set of free Schr\"oder paths of
length $2n$.
\end{theorem}

\pf Similar to that of Theorem \ref{t-c-n}.

Recall that the number of plane trees with $n$ edges and $i$
leaves is given by the Narayana number \cite{stanley}
\[ N_{n,i}=\frac{1}{n}{n
\choose i}{n \choose i-1}. \]
It follows that that the number of
leaf-colored doubly rooted plane trees equals
\begin{equation*}
\sum_{i=1}^n\left[(n+1-i)2^iN_{n,i}+i2^{i-1}N_{n,i}\right]
=\sum_{i=1}^n(2n+2-i)2^{i-1}N_{n,i}.
\end{equation*}
On the other hand, it is easy to see that the number of free
Schr\"oder paths of length $2n$ is given by the summation
\begin{equation*}
\sum_{i=0}^n{2n-i \choose i}{2n-2i \choose n-i}.
\end{equation*}
Hence Theorem \ref{t-s-n} yields the following identity:
\begin{equation}
\sum_{i=1}^n(2n+2-i)2^{i-1} N_{n,i} =\sum_{i=0}^n{2n-i \choose
i}{2n-2i \choose n-i}.
\end{equation}


By the right-to-left preorder traversal and the above
correspondence, one may determine a distinguished vertex of a
plane tree whose leaves are colored red and blue. This fact can be
restated as a Schr\"oder path analogue of the Chung-Feller theorem
 obtained by Eu, Fu and Yeh \cite{efy}.

\begin{theorem} \label{r-s}
For each Schr\"oder path $P$ from $(0,0)$ to $(2n,0)$, assign
weight $2$ to $P$ if $P$ ends with a $U$ step; otherwise $P$ is
assigned
 weight $1$. Let $m$ be a given number not exceeding $n$.
 Then the total weight of the set of free Schr\"oder paths
of length $2n$ with $m$ flaws is always the Schr\"oder number
$r_n$.
\end{theorem}

If a free Schr\"oder path ends with an up step, then the
corresponding subtree $T'$ is empty and we have that the
distinguished vertex is a leaf. There are now two possible ways to
color it, hence we assign weight $2$ to this kind of Schr\"oder
paths.

Using plane trees, we may reinterpret the above theorem as
follows.

\begin{theorem} Let $T$ a plane tree with $n$ edges. Assume that the vertices of
$T$ are labelled by $0, 1, 2, \ldots, n$ according to the
right-to-left preorder traversal. Let $w$ be a vertex labelled by
$m$. Let $T'$ be a leaf-colored doubly rooted plane tree $T$ with
$w$ being the distinguished vertex. Then by the correspondence
between leaf-colored doubly rooted plane trees and free Schr\"oder
paths, $T'$ corresponds to a free Schr\"oder path with $m$ flaws.
\end{theorem}

From the above theorem, we immediately get the following
refinement of Eu, Fu and Yeh \cite{efy}.

\begin{theorem}\label{refines}
For $0\leq k \leq m \leq n$, there is a bijection between the set
of free Schr\"oder path of length $2n$ with $m$ flaws in $k$ flaw
blocks and the set of leaf-colored doubly rooted plane trees of
$n$ edges with stem size $k$ and $m$ prefix edges.
\end{theorem}

\begin{figure}[h,t]
\begin{center}
\setlength{\unitlength}{1.3mm}
\begin{picture}(95,45)
\put(10,40){\circle*{.8}}\put(10,40){\line(-2,-1){10}}
\put(10.5,40.5){\tiny{$0$}}\put(20.5,35.5){\tiny{$1$}}
\put(9.8,40){\line(-4,-5){4}} \put(9.9,40){\line(-4,-5){4}}
\put(10,40){\line(-4,-5){4}}\put(10.1,40){\line(-4,-5){4}}
\put(10.2,40){\line(-4,-5){4}}
\put(9.8,40){\line(4,-5){4}}\put(9.9,40){\line(4,-5){4}}
\put(10,40){\line(4,-5){4}}\put(10.1,40){\line(4,-5){4}}
\put(10.2,40){\line(4,-5){4}}
\put(9.7,40){\line(2,-1){10}}\put(9.8,40){\line(2,-1){10}}
\put(9.9,40){\line(2,-1){10}}\put(10,40){\line(2,-1){10}}
\put(10.1,40){\line(2,-1){10}}\put(10.2,40){\line(2,-1){10}}
\put(10.3,40){\line(2,-1){10}}\put(10.4,40){\line(2,-1){10}}
\multiput(0,35)(6,0){2}{\circle*{.8}}
\multiput(14,35)(6,0){2}{\circle*{.8}}

\put(0,35){\line(0,-1){5}} \put(0,30){\circle*{.8}}

\put(6.8,34.5){\tiny{$5$}}
\put(5.8,35){\line(-3,-5){3}}\put(5.9,35){\line(-3,-5){3}}
\put(6,35){\line(-3,-5){3}}\put(6.1,35){\line(-3,-5){3}}
\put(6.2,35){\line(-3,-5){3}}
\put(5.8,35){\line(3,-5){3}}
\put(5.9,35){\line(3,-5){3}}\put(6,35){\line(3,-5){3}}
\put(6.1,35){\line(3,-5){3}}\put(6.2,35){\line(3,-5){3}}
\put(3,30){\circle*{.8}} \put(9,30){\circle*{.8}}

\put(14.5,35){\tiny{$2$}}\put(10.7,28){\tiny{$4$}}\put(16.7,28){\tiny{$3$}}
\put(13.8,35){\line(-3,-5){3}}\put(13.9,35){\line(-3,-5){3}}
\put(14,35){\line(-3,-5){3}}\put(14.1,35){\line(-3,-5){3}}
\put(14.2,35){\line(-3,-5){3}}
\put(13.8,35){\line(3,-5){3}}\put(13.9,35){\line(3,-5){3}}
\put(14,35){\line(3,-5){3}}\put(14.1,35){\line(3,-5){3}}
\put(14.2,35){\line(3,-5){3}}
\put(11,30){\circle*{.8}} \put(17,30){\circle*{.8}}

\put(3,30){\line(-3,-5){3}} \put(3,30){\line(3,-5){3}}
\put(0,25){\circle*{.8}} \put(6,25){\circle*{.8}}

\put(8.9,30.5){\tiny{$6$}}\put(8.5,23){\tiny{$7$}}\put(4,29.5){\tiny{$8$}}
\put(8.8,30){\line(0,-1){5}}\put(8.9,30){\line(0,-1){5}}
\put(9,30){\line(0,-1){5}}\put(9.1,30){\line(0,-1){5}}
\put(9.2,30){\line(0,-1){5}}
\put(9,25){\circle*{.8}}

\put(6,25){\line(0,-1){5}} \put(6,20){\circle*{.8}}

\put(6.7,24.5){\tiny{$9$}}\put(6.7,20){\tiny{$10$}}\put(.7,24.2){\tiny{$11$}}
\put(.8,34.2){\tiny{$12$}}\put(.4,29.4){\tiny{$13$}}

\put(-2.5,29){\scriptsize{$B$}}\put(-2.5,24){\scriptsize{$B$}}
\put(3.5,19){\scriptsize{$R$}}\put(9.5,24){\scriptsize{$B$}}
\put(11.5,29){\scriptsize{$B$}}\put(17.5,29){\scriptsize{$R$}}
\put(20.5,33.7){\scriptsize{$B$}}

\put(22,30){$\longleftrightarrow$}
\put(40,40.2){\circle{1}} \put(40,40){\line(-2,-1){10}}
\put(40,40){\line(-4,-5){4}}\put(39.8,40){\line(-4,-5){4}}\put(40.1,40){\line(-4,-5){4}}
\put(39.9,40){\line(-4,-5){4}}\put(40.2,40){\line(-4,-5){4}}
\put(40,40){\line(4,-5){4}} \put(40,40){\line(2,-1){10}}
\multiput(30,35)(6,0){2}{\circle*{.8}}
\multiput(44,35)(6,0){2}{\circle*{.8}}

\put(30,35){\line(0,-1){5}} \put(30,30){\circle*{.8}}

\put(36,35){\line(-3,-5){2.7}}\put(35.8,35){\line(-3,-5){2.7}}\put(36.1,35){\line(-3,-5){2.8}}
\put(35.9,35){\line(-3,-5){2.7}}\put(36.2,35){\line(-3,-5){2.8}}
\put(36,35){\line(3,-5){3}} \put(33,30.1){\circle{1}}
\put(39,30){\circle*{.8}}

\put(44,35){\line(-3,-5){3}} \put(44,35){\line(3,-5){3}}
\put(41,30){\circle*{.8}} \put(47,30){\circle*{.8}}

\put(33,30){\line(-3,-5){3}} \put(33,30){\line(3,-5){3}}
\put(30,25){\circle*{.8}} \put(36,25){\circle*{.8}}

\put(39,30){\line(0,-1){5}} \put(39,25){\circle*{.8}}

\put(36,25){\line(0,-1){5}} \put(36,20){\circle*{.8}}

\put(28.5,28){\scriptsize{$B$}}\put(30.5,24){\scriptsize{$B$}}
\put(36.5,19){\scriptsize{$R$}}\put(39.5,24){\scriptsize{$B$}}
\put(41.5,29){\scriptsize{$B$}}\put(47.5,29){\scriptsize{$R$}}
\put(50.5,34){\scriptsize{$B$}}

\put(52,30){$\longleftrightarrow$}


\put(61,40){\line(0,-1){10}}
\multiput(61,40)(0,-5){3}{\circle*{.8}} \put(60,23){\small{$L_1$}}
\put(61.5,29){\scriptsize{$B$}}

\put(72,40){\circle*{.8}} \put(72,40){\line(-3,-5){3}}
\put(72,40){\line(3,-5){3}} \multiput(69,35)(6,0){2}{\circle*{.8}}
\put(69,35){\line(-3,-5){3}} \put(69,35){\line(3,-5){3}}
\multiput(66,30)(6,0){2}{\circle*{.8}}\put(69,23){\small{$R_1$}}
\put(66.5,29){\scriptsize{$B$}}\put(72.5,29){\scriptsize{$R$}}
\put(75.5,34){\scriptsize{$B$}}

\put(81,35){$\emptyset$} \put(80,23){\small{$L_2$}}

\put(87,40){\line(0,-1){10}}
\multiput(87,40)(0,-5){3}{\circle*{.8}} \put(86,23){\small{$R_2$}}
\put(87.5,29){\scriptsize{$B$}}

\put(95,40){\circle*{.8}} \put(95,40){\line(-3,-5){3}}
\put(95,40){\line(3,-5){3}} \multiput(92,35)(6,0){2}{\circle*{.8}}
\put(98,35){\line(0,-1){5}}
\put(98,30){\circle*{.8}}\put(95,23){\small{$T'$}}
\put(92.5,34){\scriptsize{$B$}}\put(98.5,29){\scriptsize{$R$}}

\put(0,7.9){$\longleftrightarrow$}

\multiput(10,9)(1,0){90}{\line(1,0){0.5}}

\put(15,9){\line(1,1){3}}\put(18,12){\line(1,0){6}}
\put(24,12){\line(1,-1){9}} \put(33,3){\line(1,0){6}}
\put(39,3){\line(1,-1){3}} \put(42,0){\line(1,1){6}}
\put(48,6){\line(1,0){6}}\put(57,9){\line(1,-1){6}}
\put(63,3){\line(1,0){6}}
\put(69,3){\line(1,1){6}}\put(75,9){\line(1,0){6}}
\put(81,9){\line(1,1){6}}\put(87,15){\line(1,-1){6}}

\put(15,9){\circle*{.8}}\put(18,12){\circle*{.8}}
\put(24,12){\circle*{.8}}\put(27,9){\circle*{.8}}\put(30,6){\circle*{.8}}
\put(33,3){\circle*{.8}}\put(39,3){\circle*{.8}}
\put(42,0){\circle*{.8}}\put(45,3){\circle*{.8}}\put(48,6){\circle*{.8}}
\put(54,6){\circle*{.8}}\put(57,9){\circle*{.8}}
\put(60,6){\circle*{.8}}\put(63,3){\circle*{.8}}
\put(69,3){\circle*{.8}}\put(72,6){\circle*{.8}}\put(75,9){\circle*{.8}}
\put(81,9){\circle*{.8}}
\put(84,12){\circle*{.8}}\put(87,15){\circle*{.8}}
\put(90,12){\circle*{.8}}\put(93,9){\circle*{.8}}

\put(27.3,9){\line(1,-1){3}}\put(54.3,6){\line(1,1){3}}
\put(57.3,9){\line(1,-1){3}}\put(72.3,6){\line(1,1){3}}
\put(27.2,9){\line(1,-1){3}}\put(54.2,6){\line(1,1){3}}
\put(57.2,9){\line(1,-1){3}}\put(72.2,6){\line(1,1){3}}
\put(27.1,9){\line(1,-1){3}}\put(54.1,6){\line(1,1){3}}
\put(57.1,9){\line(1,-1){3}}\put(72.1,6){\line(1,1){3}}
\put(26.9,9){\line(1,-1){3}}\put(53.9,6){\line(1,1){3}}
\put(56.9,9){\line(1,-1){3}}\put(71.9,6){\line(1,1){3}}
\end{picture}
\end{center}
\caption{Leaf-colored plane trees and free Schr\"oder paths}
\label{sf}
\end{figure}
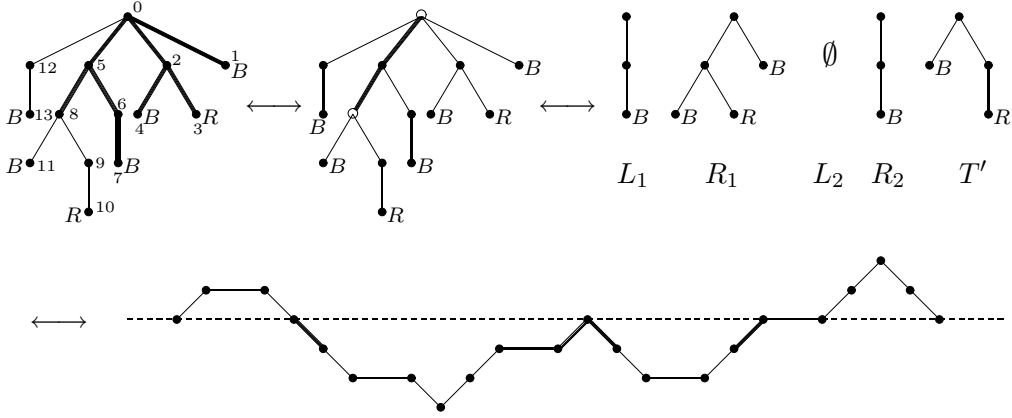

An example of the above bijection between leaf-colored doubly
rooted plane trees and free Schr\"oder paths is illustrated in
Figure~ \ref{sf}.

To conclude this section, we use the butterfly decomposition to
obtain a formula for the total weight of leaf-colored doubly
rooted plane trees of $n$ edges with stem size $k$ and $m$ prefix
edges.  Let $S$ be the generating function of the Schr\"oder
numbers as given by the equation $S=1+xS+xS^2$. Then the total
weight of leaf-colored doubly rooted plane trees of $n$ edges with
stem size $k$ and $m$ prefix edges equals
\[
2\cdot[x^m]x^kS^k\cdot[x^{n-m}]S^k+[x^m]x^kS^k\cdot[x^{n-m}]S^k(S-1)\]
which can be rewritten as $[x^{m-k}]S^k[x^{n-m}](S^{k+1}+S^k)$.
Let
\begin{equation}  \label{fnk}
 a(n,k)=[x^n]S^k.
\end{equation}
 Set $a(0,k)=1$. When $n\geq 1$, using
the Lagrange inversion formula \cite{stanley}, we obtain that
\begin{equation} \label{fk}
a(n,k)=
 \frac{k}{n}\sum_{i=0}^{n-1}2^{i+1}{n+k-1 \choose i}{n
\choose i+1}  .
\end{equation}

Note that $a(n,1)$ reduces to the Schr\"oder number $r_n$.

\begin{corollary}
For $0< k \leq m \leq n$, the total weight of free Schr\"oder
paths of length $2n$ with $m$ flaws and $k$ flaw blocks equals
\[
a(m-k,k)\cdot[a(n-m,k+1)+a(n-m,k)].
\]
\end{corollary}

\section{Two Involutions}

In this section, we present two parity reversing involutions on
free Dyck paths and free Schr\"oder paths, where the parity is
defined as the parity of the number of its flaw blocks. We also
derive two identities based on the computation via the butterfly
decomposition.

\begin{theorem} For $n\geq 1$, there is a parity reversing involution on the set of free Dyck
paths of length $2n$, which leads to the following identity
\begin{equation}\label{di}
\sum_{i=0}^n(-1)^i\frac{2i+1}{2n+1}{2n+1 \choose n-i}=0.
\end{equation}
\end{theorem}

\pf Let  $P$ be a free Dyck path of length $2n$. Let $\phi$ be the
desired involution. If $P$ ends with an elevated Dyck path, then
we construct $\phi(P)$ by reflecting the last elevated Dyck path
with respect to the $x$-axis. If $P$ ends with a negative elevated
Dyck path, then  $\phi(P)$ is obtained by reflecting the last
negative elevated Dyck path with respect to the $x$-axis. Clearly,
$\phi$ is a parity reversing involution. By the butterfly
decomposition,  the number of free Dyck paths with $i$ flaw blocks
equals the number of doubly rooted plane trees with stem size $i$,
that is, $[x^{n-i}]C^{2i+1}$. Hence the relation (\ref{di})
follows from (\ref{xn}). \qed

We also have an involution on free Schr\"oder paths.

\begin{theorem}
For $n\geq 1$, there is a parity reversing involution on the set
of free Schr\"oder paths of length $2n$ containing at least one up
step. So we have the following identity on $a(n,k)$ as defined by
(\ref{fk}):
\begin{equation}\label{si}
\sum_{i=0}^n(-1)^i\,a(n-i,2i+1)=1.
\end{equation}
\end{theorem}

\pf Let $P$ be a free Schr\"oder path which contains at least one
up step.  Let $Q$ be the last segment of $P$ which is an elevated
Schr\"oder path or a negative elevated Schr\"oder path. Note that
 $Q$ may be followed by some horizontal steps in $P$.
We reflect $Q$ with respect to the $x$-axis to get a free
Schr\"oder path. Clearly, the resulting path contains at least one
up step. It is easy to see that  this construction is reversible
and parity reversing. By the correspondence given in Theorem
\ref{refines}, the number of free Schr\"oder paths of length $2n$
with $i$ flaw blocks equals $[x^{n-i}]S^{2i+1}$, that is, $a(n-i,
2i+1)$. Therefore, identity (\ref{si}) follows from the involution
and the fact that the only Schr\"oder path not affected by the
involution is the path consisting of only horizontal steps. \qed

\section{Chains in Plane Trees}

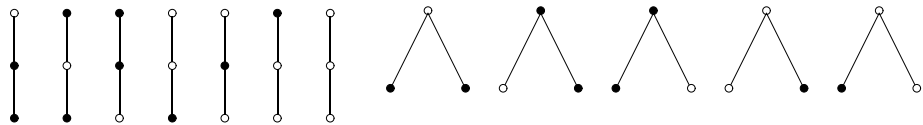
\begin{figure}[h,t]
\begin{center}
\setlength{\unitlength}{1mm}
\begin{picture}(120,20)

\put(0,0){\line(0,1){13.5}} \put(0,14){\circle{1}}
\multiput(0,0)(0,7){2}{\circle*{1}}

\put(7,0){\line(0,1){6.5}} \put(7,7.5){\line(0,1){6.5}}
\put(7,7){\circle{1}} \multiput(7,0)(0,14){2}{\circle*{1}}

\put(14,0.5){\line(0,1){13.5}} \put(14,0){\circle{1}}
\multiput(14,7)(0,7){2}{\circle*{1}}

\put(21,0){\line(0,1){6.5}} \put(21,7.5){\line(0,1){6}}
\put(21,0){\circle*{1}} \multiput(21,7)(0,7){2}{\circle{1}}

\put(28,0.5){\line(0,1){13}} \put(28,7){\circle*{1}}
\multiput(28,0)(0,14){2}{\circle{1}}

\put(35,0.5){\line(0,1){6}} \put(35,7.5){\line(0,1){6}}
\put(35,14){\circle*{1}} \multiput(35,0)(0,7){2}{\circle{1}}

\put(42,0.5){\line(0,1){6}} \put(42,7.5){\line(0,1){6}}
\multiput(42,0)(0,7){3}{\circle{1}}


\put(55,14){\line(-1,-2){5}} \put(55,14){\line(1,-2){5}}
\put(55,14.3){\circle{1}} \multiput(50,4)(10,0){2}{\circle*{1}}

\put(70,14){\line(-1,-2){4.8}} \put(70,14){\line(1,-2){5}}
\put(70,14.3){\circle*{1}} \put(65,4){\circle{1}}
\put(75,4){\circle*{1}}

\put(85,14){\line(-1,-2){5}} \put(85,14){\line(1,-2){4.8}}
\put(85,14.3){\circle*{1}} \put(80,4){\circle*{1}}
\put(90,4){\circle{1}}

\put(100,14){\line(-1,-2){4.8}} \put(100,14){\line(1,-2){5}}
\put(100,14.3){\circle{1}} \put(95,4){\circle{1}}
\put(105,4){\circle*{1}}

\put(115,14){\line(-1,-2){5}} \put(115,14){\line(1,-2){4.8}}
\put(115,14.3){\circle{1}} \put(110,4){\circle*{1}}
\put(120,4){\circle{1}}

\end{picture}
\end{center}
\caption{ Chains in plane trees with $2$ edges}
\label{chainexample}
\end{figure}

Let us recall that a {\em chain} of a plane tree is a selection of
vertices on a path from the root to a leaf.
 The \emph{size} of a chain  is defined as the number of vertices in the
chain.  Let $Q_{n}$ be the number of nonempty chains in all plane
trees with $n$ edges. A tree with $n$ edges may have as many as
$2^{n+1}-1$ non-empty chains and as few as $2n+1$. The twelve
chains in plane trees with $2$ edges are illustrated in
Figure~\ref{chainexample}, where those empty circles stand for
vertices in chains and black circles stand for normal vertices.
For instance, the last structure of Figure~\ref{chainexample} has
a chain of size $2$.

 The main result of this section is a combinatorial interpretation
 of the generating function for the number of chains in plane
 trees obtained by Klazar \cite{klazar}. We also obtain a
one-to-one correspondence between the set of chains in plane trees
with $n$ edges and the set of tricolored plane trees with $n$
edges. Klazar \cite{klazar} derived the following generating
function for the number of chains in plane trees with $n$ edges:
\begin{equation}\label{g-k}
\frac{C}{1-2xC^{2}}=1+3x+12x^{2}+51x^{3}+222x^{4}+978x^{5}+\cdots
.
\end{equation}
Note that here we use a slightly different formulation of the
generating function $C$ from that used by Klazar \cite{klazar}.

We now give a combinatorial proof of the fact that the generating
function of chains in plane trees with $n$ edges equals
  $\frac{C}{1-2xC^2}$.  Let $T$ be a plane tree and $Q$ be a chain
  of $T$. Suppose $w$ is the vertex in $Q$ such that the path
  $v_1v_2\cdots v_kw$ from the root of $T$ to $w$ contains all the vertices
  in $Q$. Moreover, we color the
  vertex $v_i$ with the white color if it belongs to $Q$; otherwise,
  we color $v_i$ with the black color.  Such a coloring scheme
  leads to the following bijection.

  \begin{theorem} There is a one-to-one correspondence between the
   set of chains in plane trees with $n$ edges and the set of
   doubly rooted plane trees in which the vertices on the path
   from the root to the distinguished vertex (but not including
   the distinguished vertex) are colored with two colors.
  \end{theorem}

Using the above theorem and the butterfly decomposition of doubly
rooted plane trees, we obtain the generating function of Klazar.

Motivated by the following relation
\begin{equation}\label{g-k-2}
{C \over 1-2xC^2} = {1 \over 1-3xC},
\end{equation}
we are led to establish the following bijection.

\begin{theorem}\label{identity}
There is a one-to-one correspondence between chains in plane trees
and tricolored plane trees.
\end{theorem}

\pf Let $T$ be a plane tree and $Q$ be a chain of $T$. Let
$v_1v_2\cdots v_kw$ be the path from the root to the vertex $w$,
where $w$ is the last vertex in the chain. Suppose that $(L_1,
R_1; L_2, R_2; \ldots; L_k,R_k;  T')$ is the butterfly
decomposition of $T$. Let $T_i$ be the planted plane tree obtained
from $R_i$ by adding a root. Coloring $L_i$ and $T'$ red, and
color $T_i$ white if the vertex $v_i$ contained in $T_i$ is a
chain vertex, otherwise color $T_i$ black. Identify their roots as
the root of the corresponding tricolored plane tree, and set the
subtrees of the root as
\[ L_1 \, T_1 \, L_2\, T_2\, \cdots \, L_k \, T_k \, T'.\]
The reverse procedure is easy to construct. This completes the
proof. \qed

\begin{figure}[h,t]
\begin{center}
\setlength{\unitlength}{1mm}
\begin{picture}(140,35)

\put(3,30){\circle*{1}}\put(3,30){\line(-3,-5){3}}
\put(3,30){\line(3,-5){2.7}}
\put(0,25){\circle*{1}}\put(6,25.3){\circle{1}}

\put(6,25){\line(-3,-5){3}}\put(3,20){\line(0,-1){5}}
\put(6,25){\line(3,-5){3}}\put(9,20){\line(0,-1){4.3}}
\multiput(3,20)(0,-5){2}{\circle*{1}} \put(9,20){\circle*{1}}
\put(9,15.3){\circle{1}}

\put(9,15){\line(-1,-1){5}}
\put(9,15){\line(0,-1){4.3}}\put(9,15){\line(1,-1){5}}
\multiput(4,10)(10,0){2}{\circle*{1}} \put(9,10.3){\circle{1}}

\put(4,10){\line(-3,-10){1.5}}\put(4,10){\line(3,-10){1.5}}
\multiput(2.5,5)(3,0){2}{\circle*{1}}

\put(9,10){\line(-3,-10){1.5}}\put(9,10){\line(3,-10){1.5}}
\multiput(7.5,5)(3,0){2}{\circle*{1}}

\put(10.5,5){\line(0,-1){5}} \put(10.5,0){\circle*{1}}

\put(14,10){\line(0,-1){5}} \put(14,5){\circle*{1}}

\put(2.5,5){\line(-3,-10){1.5}}\put(2.5,5){\line(3,-10){1.5}}
\multiput(1,0)(3,0){2}{\circle*{1}}

\put(16,15){$\longleftrightarrow$}

\put(4,30){\tiny{$v_1$}}\put(7,25){\tiny{$v_2$}}
\put(10,20){\tiny{$v_3$}}\put(10,15){\tiny{$v_4$}}
\put(10,10){\tiny{$w$}}

\put(27,25){\line(0,-1){5}} \multiput(27,25)(0,-5){2}{\circle*{1}}
\put(25,5){\small{$L_1$}}

\put(32,20){$\emptyset$} \put(31,5){\small{$R_1$}}

\put(40,25){\line(0,-1){10}}
\multiput(40,25)(0,-5){3}{\circle*{1}} \put(38,5){\small{$L_2$}}

\put(46,20){$\emptyset$} \put(45,5){\small{$R_2$}}

\put(53,20){$\emptyset$} \put(52,5){\small{$L_3$}}

\put(60,20){$\emptyset$} \put(59,5){\small{$R_3$}}

\put(73,25){\line(0,-1){5}} \multiput(73,25)(0,-5){2}{\circle*{1}}
\put(73,20){\line(-3,-5){3}}\put(73,20){\line(3,-5){3}}
\multiput(70,15)(6,0){2}{\circle*{1}}
\put(70,15){\line(-3,-5){3}}\put(70,15){\line(3,-5){3}}
\multiput(67,10)(6,0){2}{\circle*{1}} \put(70,5){\small{$L_4$}}

\put(80,25){\line(0,-1){10}}
\multiput(80,25)(0,-5){3}{\circle*{1}} \put(79,5){\small{$R_4$}}

\put(88,25){\circle*{1}} \put(88,25){\line(-3,-5){3}}
\put(88,25){\line(3,-5){3}} \multiput(85,20)(6,0){2}{\circle*{1}}
\put(91,20){\line(0,-1){5}}
\put(91,15){\circle*{1}}\put(87,5){\small{$T'$}}

\put(96,15){$\longleftrightarrow$}

\put(125,25){\circle*{1}}
\put(125,25){\line(-2,-1){16}}\put(125,25){\line(-3,-2){12}}
\put(125,25){\line(-1,-1){8}}\put(125,25){\line(-1,-2){4}}
\put(125,25){\line(0,-1){8}}\put(125,25){\line(1,-2){4}}
\put(125,25){\line(1,-1){8}}\put(125,25){\line(3,-2){12}}
\put(125,25){\line(2,-1){16}}
\multiput(109,17)(4,0){9}{\circle*{1}}

\put(109,15){\tiny{$R$}}\put(113,15){\tiny{$B$}}
\put(117,15){\tiny{$R$}}\put(121,15){\tiny{$W$}}
\put(125,15){\tiny{$B$}}\put(129.5,15){\tiny{$R$}}
\put(133.5,15){\tiny{$W$}}\put(137,15){\tiny{$R$}}
\put(141,15){\tiny{$R$}}

\put(117,17){\line(0,-1){8}}  \put(117,9){\circle*{1}}

\put(129,17){\line(-1,-4){2}}\put(129,17){\line(1,-4){2}}
\multiput(127,9)(4,0){2}{\circle*{1}}

\put(127,9){\line(-1,-4){2}}\put(127,9){\line(1,-4){2}}
\multiput(125,1)(4,0){2}{\circle*{1}}

\put(133,17){\line(0,-1){16}}
\multiput(133,9)(0,-8){2}{\circle*{1}}

\put(141,17){\line(0,-1){8}}  \put(141,9){\circle*{1}}
\end{picture}
\end{center}
\caption{Chains  and tricolored
 plane trees} \label{c3}
\end{figure}
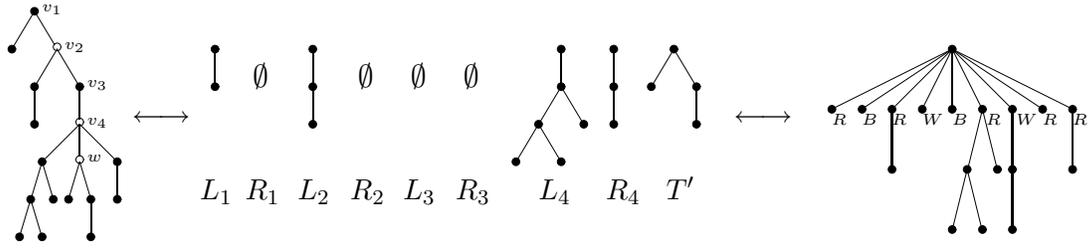

An example of the above bijection is shown in Figure~ \ref{c3}.

From the above bijection, we easily see that chains with $m$
vertices correspond to tricolored trees with $m-1$ white subtrees.
Hence as a special case of Theorem~\ref{identity}, we obtain
Theorem~\ref{doubly}.

Notice that a chain in plane trees is just a two colored path in
the butterfly decomposition. Hence we can color the vertices in
chain with $t$ colors and preserve these colors in the above
bijection. Precisely speaking, a chain is called $t$-colored if
its elements are $t$-colored. We have the following bijection.

\begin{theorem}
There is a one-to-one correspondence between the set of
$(k-2)$-colored chains in plane trees with $n$ edges and the set
of $k$-colored plane trees with $n$ edges.
\end{theorem}

The above bijection is a reflection of the following  Catalan type
identity
$$\frac{C}{1-(k-1)xC^2}=\frac{1}{1-kxC}.$$



\noindent $\mathbf{Remark.}$ By the path decomposition, the
generating function for the number of chains with $n$ edges that
end with a leaf equals
\begin{equation*}
\frac{1}{1-2xC^2}=\frac{1+\sqrt{1-4x}}{3\sqrt{1-4x}-1}
=1+2x+8x^2+34x^3+148x^4+652x^5+\cdots .
\end{equation*}
It is a new combinatorial explanation for Sequence A067336 in \cite{sloane}.

\section{Average Size of Chains }

In this section, we  use the generating function  $B$ of central
binomial coefficients as given by (\ref{db}) to study the total
size and average size of chains in plane trees with $n$ edges.
 It turns out that by a decomposition of chains we may rewrite
$\frac{C}{1-2xC^2}$ in order to give  an asymptotic formula. We
show that the average size of chains in plane trees with $n$ edges
asymptotically tends to ${ n+9 \over 6}$.

Bearing in mind that the generating function for the number of
chains of size $1$ equals the generating function $B$ of the
central binomial coefficients. We let $L^{\ast }$ be the
generating function for the number of plane trees with a
distinguished leaf. Any tree with a distinguished vertex can be
decomposed into a tree with a distinguished leaf and a subtree
rooted at the distinguished vertex. Thus we have $B=L^{\ast }C$
and $L^{\ast }=B/C$.

We now consider plane trees with at least two vertices in which
there is a distinguished leaf. Let $L$ be the generating function
of such plane trees with $n$ edges. It is easy to obtain the
following relations
 \[ L=L^{\ast
}-1=\frac{B-C}{C}=\frac{B-1}{2}.\]

\begin{property}
\label{ch3} The generating function for the total number of chains
of size $k$ in plane trees with $n$ edges equals $B\cdot \left(
\frac{B-1}{2}\right) ^{k-1}$.
\end{property}

\pf The required generating function follows from a decomposition
procedure for a plane tree with a given chain. Let $T$ be a plane
tree and $Q$ be a chain of $T$. Let $w_1, w_2, \ldots, w_k$ be the
chain vertices on the path from the root to the last vertex $w_k$.
Then $T$ can be decomposed into $k+1$ plane trees $T_1, T_2,
\ldots, T_k$, and $T'$, where $T_1$ is constructed from $T$ by
cutting off the subtrees of $w_1$, $T_2$ is obtained from the
subtree of $T$ rooted at $w_1$ by cutting off the subtrees of
$w_2$, and so on, finally $T'$ is the subtree of $T$ rooted at
$w_k$. The vertices $w_1, w_2, \ldots, w_k$ serve as distinguished
vertices in $T_1,T_2, \ldots, T_k$. The generating function for
the structure of $T_1$ equals  $L^{\ast }=B/C$, since the
distinguished vertex is allowed to coincide with the root in
$T_1$. The generating function for other $T_i$ $(2 \leq i \leq k)$
equals
 $L$ and the generating function for $T'$ equals
 $C$. Hence the required generating function equals
$L^{\ast }\cdot L^{k-1}\cdot C=B\cdot \left( \frac{B-1}{2}\right)
^{k-1}$. \qed

An interesting case arises if we look at chains of size $3$ that
include both the root and a leaf. In this case we have $L^{2}$ as
our generating function. It is easily shown that
$L^{2}=x^{2}+6x^{3}+29x^{4}+130x^{5}+\cdots $. This ubiquitous
sequence, $A008549$, also counts \cite{sloane}:
\begin{itemize}
\item The area under all Dyck paths of length $2n-2$.

\item The number of points at height one over all binomial paths
of length $2n-2$.

\item The number of inversions among all 321-avoiding permutations
in $S_{n}$.
\end{itemize}

From the above theorem, we have the following generating function
for the total number of chains in all plane trees with $n$ edges.

\begin{theorem}
The generating function for the total number of nonempty chains in
all plane trees with $ n$ edges equals $\frac{2B}{3-B}$.
\end{theorem}

\pf  We sum over $k$ to get $B\cdot \sum_{k\geq
0}\left(\frac{B-1}{2}\right)^{k}={B}\cdot{\left(1-\frac{B-1}{2}\right)^{-1}}=\frac{2B}{3-B}$.
\qed

Now we consider the asymptotic approximations. Let $H_n$ be the
total number of chains in plane trees with $n$ edges. Klazar
\cite{klazar} has shown that
\begin{equation}\label{r-n}
 H_{n}\sim \frac{1}{2}\cdot \left(
\frac{9}{2}\right)^{n}.
\end{equation}

Now we use the language of Riordan arrays
\cite{shapirotr,sprugnoli} to compute the generating function for
the total size of chains in all plane trees with $n$ edges. The
idea of Riordan arrays is represented as follows.  Given two
generating functions $g(x)=1+g_1x+g_2x^2+\cdots$ and
$f(x)=f_1x+f_2x^2+\cdots$ with $f_1\neq 0$, let
$M=(m_{i,j})_{i,j\geq 0}$ be the infinite lower triangular matrix
with nonzero entries on the main diagonal, where
$m_{i,j}=[x^i](g(x)f^j(x))$, namely, $m_{i,j}$ equals the
coefficient of $x^i$ in the expansion of the series $g(x)f^j(x)$.
If an infinite lower triangular matrix $M$ can be constructed in
this way from two generating functions $g(x)$ and $f(x)$, then it
is called a {\em Riordan array} and is denoted by
$M=(g(x),f(x))=(g,f)$.

If we multiply the matrix $M=(g,f)$ by a column vector
$(a_0,a_1,\cdots)^T$ to get a column vector $(b_0,b_1,\cdots)^T$,
then the generating functions  $A(x)$ and $R(x)$ of the sequences
$(a_0,a_1,\cdots)$ and $(b_0,b_1,\cdots)$ satisfy the following
relation
\[
R(x)=g(x)A(f(x)).\]

We now have the following generating function for the total size
of chains in plane trees with $n$ edges.

\begin{theorem}
The generating function for the total size of all chains in plane
trees with $n$ edges equals $\frac{4B}{\left( 3-B\right) ^{2}}$.
\end{theorem}

\pf Let $g(x)=B$ be the generating function for chains of size
$1$, and $f(x)=L=\frac{B-1}{2}$ be the generating function for
plane trees with at least two vertices and a distinguished leaf.
Consider the Riordan matrix $(B,L)$. The generating function of
the $j$-th $(j\geq 1)$ column is $BL^{(j-1)}$, which is the
generating function for the number of chains of size $j$. Since
the generating function of $(1,2,3,4\cdots)^T$ is
$A(x)=\frac{1}{(1-x)^2}$, it follows that the multiplication of
the Riordan matrix $(B,L)$ and the column vector
$(1,2,3,4\cdots)^T$ gives the sequence of the total size of chains
in plane trees with $n$ edges. It follows that
\[ R(x)=g(x)A(f(x))=\frac{B}{(1-L)^2}=\frac{4B}{(3-B)^2}\]
 is the
generating function for the total size of chains. The matrix
identity can be stated as
\begin{equation*}\label{rodan}
\begin{bmatrix} 1 \\2 & 1\\6 & 5 & 1\\20 & 22 & 8 & 1\\
70 & 93 & 47 & 11 & 1\\& & \cdots &&&\ddots
\end{bmatrix}\cdot
\begin{bmatrix}
1\\2\\3\\4\\5\\\vdots
\end{bmatrix}=\begin{bmatrix}
1\\4\\19\\92\\446\\\vdots
\end{bmatrix}.
\end{equation*}
This completes the proof. \qed

After some algebraic calculations we get
\[
R=\frac{\frac{5-18x}{\sqrt{1-4x}}+3}{8}\cdot \frac{1-4x}{\left(
1-\frac{9}{2}x \right) ^{2}}.
\]
Recall that Bender's lemma \cite[p.496]{bender} basically says
that if $C\left( x\right) =A\left( x\right) B\left( x\right) $ and
the radii of convergence for $A(x)$ and $B(x)$ are $\alpha $ and
$\beta $ with $\alpha $ $<$ $\beta $, then $$C_{n}\sim
A_{n}B\left( \alpha \right).$$ Let $A\left( x\right)
=\frac{1-4x}{\left( 1-\frac{9}{2}x\right) ^{2}}$ and $B\left(
x\right) =\frac{1}{8}\cdot \left(
\frac{5-18x}{\sqrt{1-4x}}+3\right)$. We have $\alpha =2/9<\beta
=1/4$ for Bender's lemma. So we have $B\left( \frac{2}{9}\right)
=\frac{3}{4}$ while $A_{n}=\frac{n+9}{2}\left( \frac{9}{2}\right)
^{n-1}$. So we obtain the following asymptotic property.

\begin{theorem} Let $R_n$ be the total size of chains in all plane trees with $n$
edges. Then we have
\begin{eqnarray} \label{s-n}
 R_n \sim \frac{n+9}{12} \left( \frac{9}{2}\right) ^{n}.
\end{eqnarray}
\end{theorem}
From Klazar's formula (\ref{r-n}) and the above formula
(\ref{s-n}) it follows that the average size of chains in planes
trees with $n$ edges approaches
\[
\lim _{n \rightarrow \infty} \frac{R_n}{H_n}=\frac{n+9}{6}.
\]
For example, \[ \frac{R_{50}}{H_{50}}=\frac{2250\,588\,
247\,788\,344\,466\,951\,528\,963\,319\,620}{ 228\,878\,
511\,199\,384\,804\,987\,952\, 173\,176\,432}\approx9.833\,1,\]
while $\frac{50+9}{6}\approx9.833\,3$.

\vskip 5mm

\noindent \textbf{Acknowledgments.} The authors thank David S.
Hough, Martin Klazar and Peter Winkler for helpful suggestions.
This work was supported by the 973 Project on Mathematical
Mechanization, the National Science Foundation, the Ministry of
Education, and the Ministry of Science and Technology of China.
The third author is partially supported by NSF grant HRD 0401697.


\end{document}